\documentclass[rqno,12pt]{amsart}

\usepackage{amsthm}
\usepackage{amsmath}
\usepackage{amssymb}
\usepackage{eucal}
\usepackage{amsfonts}
\usepackage{amsmath}
\usepackage[Lenny]{fncychap}
\usepackage{enumerate}
\usepackage{amssymb}
\usepackage[english]{babel}
\usepackage[latin1]{inputenc}
\usepackage[breaklinks=true]{hyperref}
\usepackage[all]{xy}
\usepackage{fancybox}
\usepackage{latexsym,mathrsfs,longtable,lscape,dsfont,yfonts}
\usepackage{srcltx}

\ifx\pdfpagewidth\undefined 
\usepackage[T1]{fontenc}
\else 
\usepackage[OT1]{fontenc} 
\fi 
\usepackage{amsfonts,amssymb,latexsym,longtable,amsmath}

\newtheorem{thm}{Theorem}[section]
\newcommand{\bthm}{\begin{thm}} \newcommand{\ethm}{\end{thm}}
\newtheorem{prop}[thm]{Proposition}
\newcommand{\bprp}{\begin{prop}} \newcommand{\eprp}{\end{prop}}
\newtheorem{fact}[thm]{Fact}
\newcommand{\bfct}{\begin{fact}} \newcommand{\efct}{\end{fact}}
\newtheorem{prob}[thm]{Problem}
\newcommand{\bprb}{\begin{prob}} \newcommand{\eprb}{\end{prob}}
\newtheorem{quest}[thm]{Question}
\newcommand{\bqtn}{\begin{quest}} \newcommand{\eqtn}{\end{quest}}
\newtheorem{lem}[thm]{Lemma}
\newcommand{\blem}{\begin{lem}} \newcommand{\elem}{\end{lem}}
\newtheorem{claim}[thm]{Claim}
\newcommand{\bclm}{\begin{claim}} \newcommand{\eclm}{\end{claim}}
\newtheorem{cor}[thm]{Corollary}
\newcommand{\bcor}{\begin{cor}} \newcommand{\ecor}{\end{cor}}
\newtheorem{conj}[thm]{Conjecture}
\newcommand{\bcnj}{\begin{conj}} \newcommand{\ecnj}{\end{conj}}
\theoremstyle{definition}
\newtheorem{defn}[thm]{Definition}
\newcommand{\bdfn}{\begin{defn}} \newcommand{\edfn}{\end{defn}}
\newtheorem{spec}[thm]{Specializing}
\newcommand{\bspc}{\begin{spec}} \newcommand{\espc}{\end{spec}}
\theoremstyle{remark}
\newtheorem{rem}[thm]{Remark}
\newcommand{\brem}{\begin{rem}} \newcommand{\erem}{\end{rem}}
\newtheorem{cnv}[thm]{Convention}
\newcommand{\bcnv}{\begin{cnv}} \newcommand{\ecnv}{\end{cnv}}
\newtheorem{exam}[thm]{Example}
\newcommand{\bexm}{\begin{exam}} \newcommand{\eexm}{\end{exam}}
\newcommand{\bpf}{\begin{proof}} \newcommand{\epf}{\end{proof}}

\newtheorem{thmy}{\textbf{Theorem}}

\newcommand{\Z}{\mathbb Z}

\newcommand{\N}{\mathbb N}

\renewcommand{\phi}{\varphi}
\renewcommand{\theta}{\vartheta}

\newcommand{\ga}{{\alpha}}

\newcommand{\gl}{{\lambda}}

\newcommand{\gb}{{\beta}}
\newcommand{\gga}{{\gamma}}
\newcommand{\gw}{{\omega}}
\newcommand{\gs}{{\sigma}}

\newcommand{\supp}{{\rm supp}}

\newcommand{\mkp}{\medskip}
\newcommand{\bkp}{\bigskip}

\begin{document}

\title[Homomorphic encoders of profinite abelian groups II]{Homomorphic encoders of profinite abelian groups II}

\author[M. Ferrer]{Mar\'ia V. Ferrer}
\address{Universitat Jaume I, Instituto de Matem\'aticas de Castell\'on,
Campus de Riu Sec, 12071 Castell\'{o}n, Spain.}
\email{mferrer@mat.uji.es}

\author[S. Hern\'andez]{Salvador Hern\'andez}
\address{Universitat Jaume I, INIT and Departamento de Matem\'{a}ticas,
Campus de Riu Sec, 12071 Castell\'{o}n, Spain.}
\email{hernande@mat.uji.es}

\thanks{ Research Partially supported by the Spanish Ministerio de Econom\'{i}a y Competitividad,
grant: MTM/PID2019-106529GB-I00 (AEI/FEDER, EU) and by the Universitat Jaume I, grant UJI-B2019-06}

\begin{abstract}
Let $\{G_i :i\in\N\}$ be a family of finite Abelian groups. We say that a subgroup  $G\leq \prod\limits_{i\in \N}G_i$ is
\emph{order controllable} if for every $i\in \mathbb{N}$ there is $n_i\in \mathbb{N}$ such that for each $c\in G$, there exists
$c_1\in G$ satisfying that $c_{1|[1,i]}=c_{|[1,i]}$, $supp (c_1)\subseteq [1,n_i]$, and order$(c_1)$
divides order$(c_{|[1,n_i]})$. In this paper we investigate the structure of order controllable group codes. It is
proved that if $G$ is an order controllable, shift invariant, group code over a finite abelian group $H$,
then $G$ possesses a finite canonical generating set. Furthermore, our construction also yields that $G$ is algebraically conjugate to a full group shift.

\end{abstract}

\thanks{{\em 2010 Mathematics Subject Classification.} Primary 20K25. Secondary 22C05, 20K45, 54H11, 68P30, 37B10\\
{\em Key Words and Phrases:} profinite abelian group, controllable group, order controllable group, group code,
generating set, homomorphic encoder.
}


\date{today}

\maketitle \setlength{\baselineskip}{24pt}

\section {Introduction}
In coding theory, a {\em code} refers to a set of sequences (the {\em codewords}), with good error-correcting properties,
used to transmit information over nosy channels. In communication technology most codes are linear (that is, vector spaces on a finite field)
and there are two main classes of codes: {\em block codes},
in which the codewords are finite sequences all of the same length, and {\em convolutional codes}, in which the codewords
can be infinite sequences. 
However, some very powerful codes that were first thought to be nonlinear
can be described as additive subgroups of $A^n$,
where $A$ is a cyclic abelian group (see \cite{CalHamKumSloSol,forney_trott:ie3trans93}). This fact motivated
the study of a more general class of codes.
According to Forney and Trott \cite{forney_trott:ie3trans93,forney_trott:04}, a \emph{group code} $G$ is a subgroup of
a product $$X=\prod\limits_{i\in I } G_i,$$ where each $G_i$ is a group and
the composition law is the componentwise group operation. The subgroup $$G_f:= G\cap \bigoplus _{i\in \Z} G_i$$ is called 
the \emph{finite subcode} of $G$. It may happen that all elements of $G$ has finite support, which means that $G$ coincides with $G_f$. 

If all code symbols are drawn from a common group $H$, then $G\leq H^I$
and $G$ will be called a group code over $H$ defined on $I$. 




A key point in the study of group codes is the finding of appropriate \emph{encoders} 


\bdfn
Given a group code $G$, a \emph{homomorphic encoder} is a continuous homomorphism $\Phi \colon \prod_{i\in I} H_i \to G$
that sends a full direct product of (topological) groups onto $G$.
Of special relevance are the so called \emph{noncatastrophic encoders},
that is, homomorphic continuous encoders $\alpha$ that are one-to-one and such that
$\Phi(\bigoplus_{i\in I} H_i)=G_f$ (see \cite{fagnani_zampieri:ieee96,forney_trott:ie3trans93,forney_trott:04} for some references).
\edfn


\medskip

From here on, we deal with a {\it group shift} (or {\it group code}) $G$ over a finite abelian group $H$.
That is, $G$ is a closed, \emph{shift-invariant} subgroup of the full shift group $X=H^\Z$.
Therefore, if  $\sigma\colon X\to X$ denotes the \emph{ backward shift operator}
$$\sigma[x](i):=x(i+1),\ \forall x\in X,\ i\in\Z,$$ we have that $\sigma(G)=G$.
For simplicity's sake, we denote the \emph{ forward shift operator} by $\rho$, that is
$\rho[x](i):=x(i-1),\ \forall x\in X,\ i\in\Z$. A group shift $G$ over a finite abelian group $H$ is \emph{irreducible or transitive}
if there is $x\in G$ such that the partial forward orbit $\{\sigma^n(x) : n\geq n_0\}$ is dense in $G$
for all $n_0\in\Z$.
Given two group codes $G$ and $\bar{G}$ if there is a homeomorphism (resp. topological group isomorphism)
$\Phi : G\longrightarrow \bar{G}$ so that $\gs\circ\Phi = \Phi\circ \gs$ then we say that  $G$ and $\bar{G}$ are
\emph{topologically conjugate} (resp. \emph{algebraically and topologically conjugate}) (see \cite{kitchens:1987,KitchSchmi:1989,Schmidt:1990}).

In \cite{forney:ieee1970}, Forney proved that every (linear) convolutional code is conjugate to a full shift, via a linear conjugacy.
Subsequently, it was proved by several authors (see \cite{forney_trott:ie3trans93,kitchens:1987,LoeMit:1996,MilTho:1978})
that every irreducible group shift is conjugate to a full shift. In fact, one might expect that the conjugacy was also a group homomorphism
(algebraic conjugacy). But, for group shifts, this turns to be false in general (cf. \cite{kitchens:1987,LoeMit:1996}).
In this sense, Fagnani \cite{Fagnani:1996} has obtained necessary and sufficient conditions for a group shift
to be algebraically conjugate to the full shift over a finite group. His approach is based on Pontryagin duality, which let one
reduce the question to its discrete dual group that turns out to be a finitely generated module of Laurent polynomials.

We next collect some definitions and basic facts introduced in \cite{HomEncod_I}.

\bdfn
Let $G$ be a group shift over a finite abelian group $H$. We have the following notions:
\begin{enumerate}[(1)]
\item $G$ is \emph{weakly controllable} if $G\bigcap H^{(\Z)}$ is dense in $G$; here $H^{(\Z)}$ denotes the subgroup of $H^{\Z}$
consisting of the elements  with finite support.

\item $G$ is \emph{controllable}\footnote{It is easily verified that every controllable group code $G$ is irreducible (see \cite{kitchens:1987}).}
(equivalently {\it irreducible or transitive}) if 
there is a positive integer $n_c$ such that for each $g\in G$, there exists
$g_1\in G$ such that $g_{1|(-\infty,0]}=g_{|(-\infty,0]}$ and $g_{1|]n_c,+\infty)}=0$
(we assume that $n_c$ is the least integer satisfying this property).
Remark that this property implies the existence of $g_2:=g-g_1\in G$ such that $g=g_1+g_2$,
$supp (g_1)\subseteq (-\infty,n_c]$ and $supp(g_2)\subseteq [1,+\infty[$.

\item $G$ is \emph{order controllable} if 
there is a positive integer $n_o$ such that for each $g\in G$, there exists
$g_1\in G$ such that $g_{1|(-\infty,0]}=g_{|(-\infty,0]}$, $supp (g_1)\subseteq (-\infty,n_o]$, and order$(g_{1|[1,n_0]})$
divides order$(g_{|[1,n_0]}$
(we assume that $n_o$ is the least natural number satisfying this property).
Again, this implies the existence of $g_2\in G$ such that $g=g_1+g_2$,
$supp(g_2)\subseteq [1,+\infty[$, and order$(g_2)$ divides order$(g)$.
Here, the order of $g$ is taken in the usual sense, as an element of the group $G$.
\end{enumerate}
\edfn
\mkp

We now state our main result.

\bthm\label{theorem_B}
\emph{Let $G$ be an order controllable group shift 
over a finite abelian group $H$. Then there is a noncatastrophic isomorphic encoder for $G$. As a consequence
$G$ is algebraically and topologically conjugate to a full group shift.}
\ethm

\section{Group shifts}

In this section, we apply the result accomplished in \cite[Theorem 3.2]{HomEncod_I} in order
to prove that order controllable group shifts over a finite abelian group possess canonical generating sets.
Furthermore, our construction also yields that they are algebraically conjugate to a full group shift.

In the sequel $H^{(\Z )}$ will denote the subgroup of $H^\Z$ consisting of all elements  with finite support.


\bthm\label{th_encoder2}
Let $G$ be a weakly controllable, group shift over a finite abelian $p$-group $H$.
If $G[p]$ is weakly controllable, then there is a finite generating subset
$B_0:=\{x_j: 1\leq j\leq m \}\subseteq G_{f[0,\infty)]}[p]$, 
where $x_j=p^{h_j}y_j$, $y_j\in G_{f}$, 
and each $x_j$ is selected with the maximal possible height $h_j$ in $G_f$ with
$h_j\geq h_{j+1}$, $1\leq j < m$, such that the following assertions hold true:
\begin{enumerate}
\item
There is a canonically defined $\sigma$-invariant, onto, group homomorphism
$$\Phi \colon \left(\prod\limits_{1\leq j\leq m}\Z(p^{\small{h_j+1}})\right)^{{\large \Z}}\to G.$$ 

\item (($G$ is weakly rectangular and)) $\Phi$ is a
noncatastrophic, isomorphic encoder for $G$ if
there is a finite block $[0,N]\subseteq \N$ such that the set $$\{\sigma^n[x_j]_{|[0, N]}\not= 0: n\in\Z, 1\leq j\leq m \}$$
is linearly independent.

\end{enumerate}
\ethm
\bpf
$(1)$ Using that $G$ and $G[p]$ are weakly controllable, we can proceed
as in \cite[Theorem 3.2]{HomEncod_I} in order to define a subset $B_0:=\{x_1,\dots ,x_m\}\subseteq G_{f}[p]_{[0,\infty)}$
such that $\pi_{[0]}(B_0)$ forms a basis of $\pi_{[0]}(G_{[0,\infty)}[p])$ and
for each $x_j\in B_0$ there is a nonnegative integer $h_j$ and an element $y_j\in G_f$
such that $x_j=p^{h_j}y_j$, where each $x_j$ has the maximal possible height $h_j$ in $G_f$ and $h_1\geq h_2\geq \dots\geq h_m$.
Now define $$\phi_0\colon \Z(p)^m\to G[p]$$ by $$\phi_0[(\gl_1,\dots ,\gl_m)]=\gl_1x_1+\dots +\gl_mx_m$$ and, for each $n\in\Z$, $n>0$, set
$B_n:=\rho^n(B_0)\subseteq G_{f}[p]_{[n,\infty)}$ and define  
$$\phi_n\colon Z(p)^m\to G[p]$$ by $$\phi_n[(\gl_1,\dots ,\gl_m)]=\gl_1\rho^n(x_1)+\dots +\gl_m\rho^n(x_m).$$ Now, we can define
$$\oplus_{n}\,\phi_n \colon\bigoplus\limits_{n\geq 0} (\Z(p)^m)_n\to G_{f}[p]_{[0,\infty)}$$
by $$\oplus_n\,\phi_n[\sum\limits_{n\geq 0} (\gl_{1n},\gl_{2n},\dots ,\gl_{mn})]:=\sum\limits_{n\geq 0} \phi_{n}[(\gl_{1n},\gl_{2n},\dots ,\gl_{mn})],$$
where $(\Z(p)^m)_n=\Z(p)^m$ for all $n\geq 0$.

Remark that all the maps set above are well defined group homomorphisms since each of these maps involves finite sums in its definition.
Furthermore, since the range of $\phi_n$ is contained in $G_{f}[p]_{[n,\infty)}$ for all $n\geq 0$, it follows that
the map $\oplus_{n}\,\phi_n$ is continuous when its domain (and its range) is equipped with the product topology. Therefore, there is
a canonical extension of $\oplus_{n}\,\phi_n$ to a continuous group homomorphism
$$\Phi_0: \prod\limits_{n\geq 0}(\Z(p)^m)_n\to G[p]_{[0,\infty)}.$$

Now, repeating the same arguments as in \cite[Theorem 3.2]{HomEncod_I}, it follows that
$$G_f[p]_{[0,\infty)}\subseteq \Phi_0(\prod\limits_{n\geq 0}(\Z(p)^m)_n),$$
which implies that $\Phi_0$ is a continuous onto group homomorphism because
$G_f[p]_{[0,\infty)}$ is dense in $G[p]_{[0,\infty)}.$
Furthermore, using the $\sigma$-invariance of $G$,
we can extend  $\Phi_0$ canonically to continuous onto group homomorphism
$$\Phi_N: \prod\limits_{n\geq -N} (\Z(p)^m)_n\longrightarrow G[p]_{[-N,\infty)}$$
by
$$\Phi_N[\sum\limits_{n\geq -N} (\gl_{1n},\gl_{2n},\dots ,\gl_{mn})]:=\sigma^{N}[\Phi_0[\rho^N[\sum\limits_{n\geq -N} (\gl_{1n},\gl_{2n},\dots ,\gl_{mn})]]],$$
for every $N>0$. Now, if we identify $\prod\limits_{n\geq -N} (\Z(p)^m)_n$ with the subgroup
$(\prod\limits_{n\in \Z}(\Z(p)^m)_n)_{[-N,+\infty)}$, remark that $\Phi_{(N+1)}$ restricted to
$\prod\limits_{n\geq -N} (\Z(p)^m)_n$ is equal to $\Phi_N$. Therefore, we have defined a map
$$\Phi_{\infty}\colon \bigcup\limits_{N>0}\, \prod\limits_{n\geq -N} (\Z(p)^m)_n\to G[p].$$
Again, because $\bigcup\limits_{N>0}\, \prod\limits_{n\geq -N} (\Z(p)^m)_n$
is dense in $(\Z(p)^m)^\Z,$ 
it follows that we can extend $\Phi_\infty$ to a continuous onto group homomorphism
$$\Phi:(\Z(p)^m)^\Z\longrightarrow G[p].$$

Now, taking into account that $\lim\limits_{n\rightarrow \pm\infty} \sigma^n(y_j)=0$ for all $1\leq j\leq m$,
we proceed as in \cite[Theorem 3.2]{HomEncod_I} in order to lift $\Phi$ to a continuous onto group homomorphism
$$\Phi \colon \left(\prod\limits_{1\leq j\leq m}\Z(p^{h_j+1})\right)^\Z\to G.$$
This completes the proof of (1).

(2) First, we remark that repeating the proof accomplished in \cite[Theorem 3.2]{HomEncod_I}, it follows that
the sets  $\{\sigma^n[x_j] : n\in\Z, 1\leq j\leq m \}$
and  $\{\sigma^n[y_j] : n\in\Z, 1\leq j\leq m \}$ are both (linearly) independent.

Furthermore, since all elements $x_j$ ($1\leq j\leq m$) have finite support,
it follows that the set $\{\sigma^n[x_j]_{|[0, N]}\not= 0: n\in\Z, 1\leq j\leq m \}$ is finite.
Thus, using the $\sigma$-invarince of $G$, we proceed as in \cite[Theorem 3.2]{HomEncod_I} to obtain
that $\Phi$ is one-to-one.

In order to prove that $\Phi$ is noncatastrophic, that is
$\Phi[(\prod\limits_{1\leq j\leq m}\Z(p^{h_j+1})^{(\Z)}]\subseteq G_f$, first notice that $\Phi^{-1}$ is continuous,
being the inverse map a continuous one-to-one group homomorphism. Now, reasoning by contradiction,
suppose there is $w\in G_f$ such that $(\vec{\gl_n})=\Phi^{-1}(w)$ is an infinite sequence,
let us say, without loss of generality, an infinite sequence on the right side. Then we have that the sequence
$(\sigma^n(w))_{n>0}$ converges to $0$ in $G$. However, since $(\vec{\gl_n})$ is infinite on the right side,
it follows that the sequence $(\Phi^{-1}(\sigma^n(w)))_{n>0}=(\sigma^n[(\vec{\gl_n})])_{n>0}$ does not converge to $0$
in $(\prod\limits_{1\leq j\leq m}\Z(p^{(h_j+1)}))^\Z$. This contradiction completes the proof.
\epf
\medskip

\bdfn
In the sequel, a set $\{y_1,\dots ,y_{m}\}$ (resp. $\{x_1,\dots ,x_{m}\}$) that satisfies the properties
established in Theorem \ref{th_encoder2} is called \emph{topological generating set} of $G$ (resp. $G[p]$).
\edfn

Next, we are going to use the preceding results in order to characterize the existence of noncatastrophic, isomorphic encoders.
As a consequence, we also characterize when a group shift is algebraically conjugate to a full group shift.
First we need the following notions.

\bdfn
A group shift $G\subseteq X=H^\Z$ is a shift of  {\it finite type} (equivalently, is an {\it observable} group code)
if it is defined by forbidding the appearance a finite list of (finite) blocks. As a consequence, there is
$N\in\N$ such that if $x_1,x_2$ belong to $G$ and they coincide on an $N$-block $[k,\dots , k+N]$,
then there is $x\in G$ such that $x_{|(-\infty , k+N]}=x_{1|(-\infty , k+N]}$ and
$x_{|[k,\infty)}=x_{2|[k,\infty)}$. It is known that if $G$ is an irreducible group shift over a finite group $H$,
then $G$ is also a group shift of finite type (see \cite[Prop. 4]{kitchens:1987}). Moreover, since every order controllable
group shift $G$ is irreducible, it follows that order controllable group shifts are of finite type.

Given an element $x\in G_f$ with  $\supp (x)=\{i\in\Z : x(i)\not=0\},$
the first index (resp. last index) $i\in\supp (x)$ is denoted by $i_f(x)$ (resp. $i_l(x)$).
The lenth of $\supp (x)$ is defined as $|\supp (x)|:= i_l(x)-i_f(x)+1$.
\edfn

\bprp\label{Pr_p_full}
Let $G$ be a weakly controllable, group shift of finite type over a finite abelian $p$-group $H$. If $exp(H)=p$,
then there is a noncatastrophic isomorphic encoder for $G$. As a consequence
$G$ is algebraically and topologically conjugate to a full group shift.
\eprp
\bpf
First, remark that $G=G[p]$ in this case. By Theorem \ref{th_encoder2}, there is a topological generating subset
$\mathcal{B}_0:=\{x_j: 1\leq j\leq m \}\subseteq G_{f[0,\infty]}[p]=G_{f[0,\infty]}$
such that $\pi_{[0]}(\mathcal{B}_0)$ forms a basis of $\pi_{[0]}(G_{[0,\infty)})$ and
there is a canonically defined $\sigma$-invariant, onto, group homomorphism
$$\Phi \colon (\Z(p)^{\, m})^\Z\to G.$$
Furthermore, we select each element $x_j$ with minimal support in $G_{f[0,\infty)}$ and
such that $|\supp (x_1)|\leq\dots \leq |\supp (x_m)|$.

By Theorem \ref{th_encoder2} (2), it will suffice to verify that there is
a finite block $[0,N]\subseteq \N$ such that the set
$\{\sigma^n[x_j]_{|[0, N]}\not= 0: n\in\Z, 1\leq j\leq m \}$
is linearly independent.  Indeed, let $N$ be a natural number such that
$\supp (x_j)\subseteq [0,N]$ for all $1\leq j\leq m$ and satisfying the condition of being
a group shift of finite type for $G$. That is,
if $\gw_1,\gw_2$ belong to $G$ and they coincide on any $N$-block $[k,\dots , k+N]$,
then there is $w\in G$ such that $w_{|(-\infty , k+N]}=w_{1|(-\infty , k+N]}$ and
$w_{|[k,\infty)}=w_{2|[k,\infty)}$.

Reasoning by contradiction, let us suppose that there is a linear combination
$$\sum \gl_{nj}\gs^n(x_j)_{|[0,N]}=0.$$

Since the set $\{\sigma^n[x_j] : n\in\Z, 1\leq j\leq m \}$ is linearly independent, there must be
an element $u=\sigma^{n_1}[x_{j_1}]$ (for some $n_1$ and $j_1$) such that
$$\supp (u)\cap (-\infty,0)\not=\emptyset.$$

As a consequence, there exist $\{\ga_{nj}\}\subseteq \Z(p)$ such that
$$u_{|[0,N]}= \sum\limits_{(n\not=n_1,j\not=j_1)} \ga_{nj}\gs^n(x_j)_{|[0,N]}.$$

We select $u$ such that $i_f(u)$ is minimal among the elements satisfying this property.
Set $$v:=  \sum\limits_{(n\not=n_1,j\not=j_1)} \ga_{nj}\gs^n(x_j).$$
We have that $$(u-v)_{|[0,N]}=0.$$ Since $G$ is of finite type for $N$-blocks,
there exists $w\in G$ such that
$$w_{|(-\infty , N]}=(u-v)_{|(-\infty , N]}\ \hbox{and}\ w_{|[0,\infty)}=0.$$

We have that $i_f(u)\leq i_f(w)$ and $i_l(w)<i_l(u)$. Therefore, we have found an element $w\in G_f$
with $|\supp (w)|<|\supp (u)|$. Therefore, we can replace $x_{j_1}$ by $\widetilde{x}_{j_1}:=\gs^{-n_1}(w)$
 and $|\supp (\widetilde{x}_{j_1})|<|\supp (x_{j_1})|$. This is a contradiction with our previous selection
of the (ordered) set $\{x_j: 1\leq j\leq m \}$,
which completes the proof.
\epf
\bkp

\blem\label{Le_finito}
Let $G$ be an order controllable group shift over a finite abelian $p$-group $H$.
Then $G[p]$ and $p^rG$ are order controllable group shifts for all $r$ with $p^r < exp(H)$.
As a consequence, it holds that $(p^rG)_f=p^rG_f$ for all $r$ with $p^r < exp(H)$.
\elem
\bpf
It is obvious that $G[p]$ is order controllable.
Regarding the group $p^rG$,
take an arbitrary element $x=p^ry\in p^rG$. By the order controllability of $G$, there is $z\in G$ and $n_0\in\N$
such that $y_{|(-\infty,0]}=z_{|(-\infty,0]}$, $\supp(z)\subseteq (-\infty, n_0]$
and order$(z_{|[1,n_0]})$ divides order$(y_{|[1, n_0]})$. Then $p^rz\in p^rG$,
$x_{|(-\infty,0]}=p^rz_{|(-\infty,0]}$, $\supp(p^rz)\subseteq (-\infty, n_0]$
and order$(p^rz_{|[1,n_0]})$ divides order$(x_{|[1, n_0]})$.

Finally, it is clear that $p^rG_f\subseteq (p^rG)_f$. Next we check the reverse implication.

Let $y\in G$ such that $x=p^ry\in (p^rG)_f$. Then there are two integers $m, M$ such that $x\in G_{[m,M]}$.
Assume that $M\geq 0$ without loss of generality.
By order controllability, there is $z\in G$ such that $\sigma^M(y)_{|(-\infty,0]}=z_{|(-\infty,0]}$, $\supp(z)\subseteq (-\infty, n_0]$
and order$(z_{|[1,n_0]})$ divides order$(\sigma^M(y)_{|[1, n_0]})$. Hence, if $v=\sigma^{-M}(z)$, we have
$y_{|(-\infty,M]}=v_{|(-\infty,M]}$, $\supp(v)\subseteq (-\infty, M+n_0]$
and order$(v_{|[M+1,M+n_0]})$ divides order$(y_{|[M+1,M+n_0]})$. Therefore $x=p^rv$ with $v\in G_{(-\infty,M+n_o]}$.

If $m-n_o>0$, by order controllability, there is $u\in G$ such that
$v_{|(-\infty,0]}=u_{|(-\infty,0]}$, $\supp(u)\subseteq (-\infty, n_0]\subseteq (-\infty, m-1]$
and order$(u_{|[1,n_0]})$ divides order$(v_{|[1, n_0]})$. Set $w=v-u$. We have that $w\in G_{[1,M+n_o]}$
and $x=p^rw$, which yields $x\in p^rG_f$.

If $m-n_o\leq 0$, set $N=m-n_o-1$.
By order controllability, there is $u_1\in G$ such that
$\sigma^{N}(v)_{|(-\infty,0]}=u_{1|(-\infty,0]}$, $\supp(u_1)\subseteq (-\infty, n_0]$
and order$(u_{1|[1,n_0]})$ divides order$(\sigma^{N}(v)_{|[1, n_0]})$. Hence, if $u_2=\sigma^{-N}(u_1)$, we have
$v_{|(-\infty,N]}=u_{2|(-\infty,N]}$, $\supp(u_2)\subseteq (-\infty, N+n_0]\subseteq (-\infty, m-1]$
and order$(u_{2|[N+1,N+n_0]})$ divides order$(v_{|[N+1,N+n_0]})$. Set $w=v-u_2$. We have that $w\in G_{[N+1,M+n_o]}$
and $x=p^rw$, which agin yields $x\in p^rG_f$. This completes the proof.

\epf
\mkp

Let $G$ be a group shift over a finite abelian $p$-group $H$ and 
let $G/pG$ denote the quotient group defined by the map $\pi\colon G\to G/pG$. We define the subgroup
$$(G/pG)_f:= \{\pi(u) : u\in G\ \hbox{and}\ u(n)\in pH\ \hbox{for all but finitely many}\ n\in\Z .\}$$

\blem\label{Le_base}
Let $G$ be an order controllable group shift  over a finite abelian $p$-group $H$ and
let $\{x_1,\dots ,x_m\}\subseteq (pG_{f})_{[0,\infty)}$ be a topological generating set of $pG$,
where $x_i=py_i$, $y_i\in G_f$, $1\leq i\leq m$.
If $u\in G_{f}$ 
then there exist $v\in G_{f}[p]$ 
and $w \in \langle\{\gs^n(y_j) : n\in\Z,\ 1\leq j\leq m \}\rangle$ such that $u=v+w$.
\elem
\bpf
Since  $\{x_1,\dots ,x_m\}$ is a topological generating set of $pG$,
we have
$$pu = \sum\limits_{n\in\Z}\sum\limits_{i=1}^m\gl_{in}\gs^n(x_i)= \sum\limits_{n\in\Z}\sum\limits_{i=1}^m\gl_{in}p\gs^n(y_i)=
p\sum\limits_{n\in\Z}\sum\limits_{i=1}^m\gl_{in}\gs^n(y_i).$$
Furthermore, since the group shift $pG$ is of finite type and $(pG)_f=p(G_f)$ by Lemma \ref{Le_finito},
we can apply Proposition \ref{Pr_p_full} to the group shift $pG$,
in order to obtain that the sum in the equality above only involves non-null terms for
a finite subset of indices $F\subseteq \Z$. 
Therefore
$$pu = p\sum\limits_{n\in F}\sum\limits_{i=1}^m\gl_{in}\gs^n(y_i).$$
Set $$w:=  \sum\limits_{n\in F}\sum\limits_{i=1}^m\gl_{in}\gs^n(y_i)\in G_f.$$
Then $$u=w+(u-w),$$
where $w \in \langle\{\gs^n(y_j) : n\in\Z,\ 1\leq j\leq m \}\rangle$ and $p(u-w)=0.$
It now suffices to take $v:=u-w$.
\epf

\bthm\label{Th_p_full}
Let $G$ be an order controllable group shift (therefore, of finite type) over a
finite abelian $p$-group $H$. Then there is a noncatastrophic isomorphic encoder for $G$. As a consequence
$G$ is algebraically and topologically conjugate to a full group shift.
\ethm
\bpf
Using induction on the exponent of $G$, we will prove that there is topological generating set $B_0$ of $G[p]$,
where $B_0:=\{x_1,\dots ,x_m\}\subseteq (pG_{f}[p])_{[0,\infty)}$ such that
$\pi_{[0]}(B_0)$ forms a basis of $\pi_{[0]}(\,(pG[p])_{[0,\infty)})$ and
for each $x_j\in B_0$ there is an element $y_j\in G_{f}$
such that $x_j=p^{h_j}y_j$. Furthermore $G$ is algebraically conjugate to the full group shift generated by
$\Z(p^{h_1})\times\dots \Z({p^{h_m}}).$

The case $exp(G)=p$ has already be done in Proposition \ref{Pr_p_full}.
Now, suppose that the proof has been accomplished if $exp(G)=p^h$ and let us
verify it for $exp(G)=p^{h+1}$.
We proceed as follows:

First, take the closed, shift invariant, subgroup $pG$. We have that $exp(pG)=p^h$
and by the induction hypothesis, there is topological generating set $B_0$ of $pG[p]$,
where $B_0:=\{x_1,\dots ,x_m\}\subseteq (pG_{f}[p])_{[0,\infty)}$ such that
$\pi_{[0]}(B_0)$ forms a basis of $\pi_{[0]}(\,(pG[p])_{[0,\infty)})$ and
for each $x_j\in B_0$ there is an element $y_j\in pG_{f}$
such that $x_j=p^{h_j}y_j$.

Since $y_j\in pG_{f}$, there is $z_j\in G_f$ such that $y_j=pz_j$, $1\leq j\leq m$.
Furthermore, we may assume that there is a finite block $[0,N_1]\subseteq \N$ such that the set
$\{\sigma^n[y_j]_{|[0, N_1]}\not= 0: n\in\Z, 1\leq j\leq m \}$
is linearly independent. As a consequence, using similar arguments as in \cite[Theorem 3.2]{HomEncod_I},
it follows that the set $\{\sigma^n[z_j]_{|[0, N_1]}\not= 0: n\in\Z, 1\leq j\leq m \}$
also is linearly independent. Therefore
there is a canonically defined $\sigma$-invariant, onto, group homomorphism
$$\Phi \colon \left(\prod\limits_{1\leq j\leq m} \Z({p^{h_j}})\right)^\Z\to pG.$$

Now, we complete the set $B_0:=\{x_1,\dots ,x_m\}\subseteq (pG)[p]_{f[0,\infty)}$
with a finite set $B_1:=\{u_1,\dots , u_k\}\subseteq G[p]_{f[0,\infty)}$ such that
$\pi_{[0]}(B_0\cup B_1)$ is a basis
of $\pi_{[0]}(G[p])$. Remark that we must have $h(u_i)=0$ for all $1\leq i\leq k$, since
$\pi_{[0]}(B_0)$ forms a basis of $\pi_{[0]}(\,(pG[p])$. Furthermore, arguing as
in Proposition \ref{Pr_p_full}, we may assume that there is a finite block $[0,N_2]\subseteq \N$
such that the set $$E:=\{\sigma^n[u_i]_{|[0, N_2]} : \sigma^n[u_i]_{|[0, N_2]}\not= 0: n\in\Z, 1\leq i\leq k \}$$
is an independent subset of $G[p]_{|[0, N_2]}$.

Now, consider the quotient group homomorphism $$q\colon G \to G/pG$$
and remark that $G/pG$ is a group shift over $(H/pH)^\Z$. Making use of this quotient map,
we select a basis $$V_1:=\{v_1,\dots , v_k\}\subseteq G_f[p]_{[0,+\infty)}$$
satisfying the following properties:

\begin{enumerate}
\item $V_{1|[0,N_2]}\subseteq\langle\{\sigma^n[u_i]_{|[0, N_2]} :
\sigma^n[v_i]_{|[0, N_2]}\not= 0: n\in\Z, 1\leq i\leq k \}\rangle.$

\item $\pi_{[0]}(B_0\cup V_1)$ is a basis of $\pi_{[0]}(G[p])$.

\item The set $$\{\sigma^n[v_i]_{|[0, N_2]} : \sigma^n[v_i]_{|[0, N_2]}\not= 0: n\in\Z, 1\leq i\leq k \}$$
is independent.

\item Each $q(v_i)$ has the minimal possible support in $(G/pG)_f$. That is
$$|\supp (q(v_1))|\leq\dots \leq |\supp (q(v_k))|$$ where, if $\supp (q(v_i))=\{\dots , l_1,\dots ,l_{p_i}\}$,
then $|\supp (q(v_i))|:= l_{p_i}-l_1+1$.
\end{enumerate}
\mkp

It is straightforward to verify that $q(G_f)\subseteq (G/pG)_f$ and, as a consequence,
it follows that the group $G/pG$ is controllable and its controllability index is less than or equal to
the controllability index of $G$. As in Theorem \ref{th_encoder2}, the topological generating set
$\{v_1,\dots ,v_k\}\cup \{z_1,\dots ,z_m\}$ defines a continuous onto group homomorphism
$$\Phi: \left(\Z(p)^k\times\prod\limits_{1\leq j\leq m} (\Z_{p^{h_m+1}}\right)^\Z\longrightarrow G$$


By Theorem \ref{th_encoder2}, in order to proof that $\Phi$ is one-to-one, it will suffice to find
some block $[0, N]\in\Z$ such that
$$S:=\left(\{\sigma^s[v_i]_{|[0, N]}\not= 0: s\in\Z, 1\leq i\leq k \}\cup \{\sigma^n[z_j]_{|[0, N_1]}\not= 0:
n\in\Z, 1\leq j\leq m \}\right)_{|[0,N]}$$
forms and independent subset of $G_{|[0,N]}$.

Since that this property holds separately for $\{z_1,\dots ,z_m\}$ on the block $[0,N_1]$
and $\{v_1,\dots ,v_k\}$ on the block $[0,N_2]$, it will suffice to verify that if we denote by
$Y$ to the group shift generated by $\{z_1,\dots , z_m\}$ and by $U$ to
the group shift generated by $\{v_1,\dots , v_k\},$ then
there is an block $[0,N]\subseteq \Z$ such that
$$(Y\cap U)_{|[0,N]}=\{0\}.$$

Since this implies that $S_{|[0,N]}$ is an independet subset.

Indeed, take $N\geq \max (2N_1, 2N_2)$. Then, reasoning by contradiction, assume we have a sum
$$(\sum \ga_{in}\gs^n(v_i) + \sum \gb_{js}\gs^s(z_j))_{|[0,N]}=\{0\}.$$
Remark that we may assume that this sum is finite without loss of generality since $G$ is order controllable.
Then
$$p(\sum  \ga_{in}\gs^n(v_i) + \sum \gb_{js}\gs^s(z_j))_{|[0,N]}=
(\sum p\ga_{in}\gs^n(v_i) + \sum p\gb_{js}\gs^s(z_j))_{|[0,N]}=\{0\}$$
this yields
$$ \sum p\gb_{js}\gs^s(z_j)_{|[0,N]}= \sum \gb_{js}\gs^s(y_j)_{|[0,N]}=\{0\}.$$
Since $N\geq N_1$, this implies that
$$ \sum\gb_{js}\gs^s(y_j)=\{0\}.$$

\noindent This mans that $\gb_{js}=p\gga_{js}$ for every index $js$. Thus we have
$$(\sum \ga_{in}\gs^n(v_i) + \sum p\gga_{js}\gs^s(z_j))_{|[0,N]}=\{0\}.$$

\noindent Now, we select an element $\gs^n(v_i)$ such that 
$i_f(q(\gs^n(u_i)))$ is minimal among the elements satisfying this property.
Suppose without loss of generality that $\gs^n(v_i)=\gs^{n_1}v_1$ for simplicity's sake. Solving for $\gs^{n_1}v_1$ in the
equality above, we have
$$\gs^{n_1}v_{1|[0,N]}= (\sum\limits_{n\not=n_1,i\not=1} \ga'_{in}\gs^n(v_i) + \sum p\gga'_{js}\gs^s(z_j))_{|[0,N]}=\{0\}.$$

\noindent Set $$w:= \sum\limits_{n\not=n_1,i\not=1} \ga'_{in}\gs^n(v_i)$$
and set $$w_1:=\gs^{n_1}v_1-w.$$

\noindent Remark that $pw_1=0$, that is $w_1\in G[p]$ and
$$w_{1|[0,N]}= \sum p\gga'_{js}\gs^s(y_j)_{|[0,N]}\in pH.$$
Therefore $$\supp (q(w_1))\cap [0,N]=\emptyset.$$
Since $G$ is a group shift of finite type, there is $w_2\in G$ such that
$$w_{2|(-\infty, N]}= w_{1|(-\infty, N]}\ \ \hbox{and}\ \ w_{2|[0,+\infty)}= \sum p\gga'_{js}\gs^s(y_j)_{|[0,+\infty)}.$$
From the way $w_2$ has been defined, we have that $\gs^{-n_1}(w_2)\in G_f[p]_{[0,+\infty)}$ satisfies that
$$\gs^{-n_1}(w_2)_{|[0,N_2]}\in\langle\{\sigma^n[v_i]_{|[0, N_2]} :
\sigma^n[v_i]_{|[0, N_2]}\not= 0: n\in\Z, 1\leq i\leq k \}\rangle$$
and
$$|\supp (q(w_2))|\leq |\supp (q(\gs^{n_1}(v_1)))|.$$
This is a contradiction and completes the proof.
\epf
\mkp

We can now proof Theorem \ref{theorem_B}.

\bpf[\textbf{Proof of Theorem \ref{theorem_B}}]
Since every finite abelian group is the direct sum of all its nontrivial p-subgroups,
the proof follows from Theorem \ref{Th_p_full}, in like manner as \cite[Theorem A]{HomEncod_I}
follows from \cite[Theorem 3.2]{HomEncod_I}.
\epf
\mkp

\noindent {\bf QUESTION:} Under what conditions is it possible to extend Theorem \ref{theorem_B} to non-abelian groups?
\mkp




\end{document}